\documentclass[a4paper,12pt]{article}
\usepackage{moreverb}
\usepackage{lineno}
\usepackage{graphicx}
\usepackage[english]{babel}
\usepackage[latin1]{inputenc}
\usepackage{amsmath}
\usepackage{amssymb}
\usepackage{multirow}
\usepackage{amsthm}
\usepackage{amsfonts}
\usepackage{graphicx}
\theoremstyle{plain}

\newtheorem{definition}{Definition}[section]

\title{Non Asymptotic Performance of Some Markov Chain Order Estimators}
\author{ A.R. Baigorri\thanks{baig@unb.br}\\ Mathematics Department \\ UnB }

\date{October 28, 2012}

\begin{document}
\maketitle
\linenumbers

\begin{abstract}
In what follows we study  non asymptotic behavior of different  well known  estimators \textbf{\textit{AIC}}(\cite{Tong}), 
\textbf{\textit{BIC}}(\cite{Schwarz}) and \textbf{\textit{EDC}}(\cite{Zhao,Dorea}) in contrast with the Markov chain order estimator, named as 
Global Depency Level-\textbf{\textit{GDL}}(\cite{Baigorri}).

The estimator \textit{GDL}, is based on a  different principle which makes it behave 
in a quite different form. It is strongly consistent and  more efficient than \textit{AIC}(inconsistent), outperforming 
the well established and consistent \textit{BIC} and \textit{EDC}, mainly on relatively small samples.

The estimators  mentioned above mainly consist in the evaluation of the Markov chain's sample by different multivariate \textit{deterministic} functions. 
The log likelihood approach, as in (\ref{log_likelihood_difference}),

\begin{linenomath*}
\begin{equation} \label{log_likelihood_difference}
\hspace{-15pt}\boxed{\,\,\,
\mathbb{L}_{[(n \,,\, a_1^{\kappa})\,(i,k)]}\Big(\mathbb{X}_{a_1^\kappa}(i,k)\Big)
= {\cal L}_{1} \Big( {\cal L}_{2} -  \pi(i) \,\, \mathbb{X}_{a_1^\kappa}(i,k) \, \log \mathbb{X}_{a_1^\kappa}(i,k) \Big )
\,\,\,}
\end{equation}
\end{linenomath*}

with deterministic function
\begin{linenomath*}
\begin{equation*}
\hspace{-15pt}\boxed{\,\,\,
\textit{\textbf{L}}_{[(n \,,\, a_1^{\kappa})]} 
= {\cal L}_{1} \Big( {\cal L}_{2} - \pi(i) \, x(i,j) \, \log x(i,j),\,\,\,\, {\cal L}_{1}  = const., {\cal L}_{2}  = const. 
\,\,\,}
\end{equation*}
\end{linenomath*}
\vspace{10pt}

or, the GDL approach, as in (\ref{function_gdl_for_taylor})),
\begin{linenomath*}
\begin{equation*}
\boxed{\,\,
\mathbb{G}_{[({n,\,a_1^\kappa})\,(i,j)]}\big(...,\mathbb{X}_{a_1^\kappa}(s,t),...\big) =  \sum_{i=1}^{m} \sum_{j=1}^{m}
\textit{\textbf{G}}_{[(n \,,\, a_1^{\kappa})\,(i,j)]} \big(...,\mathbb{X}_{a_1^\kappa}(s,t),...\big) 
\,\,}
\end{equation*}
\end{linenomath*}

with deterministic function
\begin{linenomath*}
\begin{equation*}
\boxed{\,\,
\textit{\textbf{G}}_{[(n \,,\, a_1^{\kappa})\,(i,j)]} \big(...,x(s,t),...\big) = 
\frac{ \Big( x(i,j) - 
       \big[ \sum_{t=1}^m  x(i,t) \big] \big[ \sum_{s=1}^m  x(s,j) \big] \Big)^2}
     {( \sum_{t=1}^m  x(i,t) ) ( \sum_{s=1}^m  x(s,j) )}. \,\,}
\end{equation*}
\end{linenomath*}
shall be analized in \textit{\textbf{Section \ref{relevant_functions}}},
exhibiting  different structural  properties. 
It will become clear the intimate differences existing between the variance of both estimators, which induce 
quite dissimilar performance, mainly for samples of moderated sizes.
\end{abstract}

\vspace{20pt}
\section{Introduction}
A Markov Chain is a discrete stochastic process $\mathbb{X}= {\{ X_n \}}_{n \geq 0} $ with   state space $E$,
cardinality $\vert E \vert < \infty$ for which there is a $ k \geq 1 $ such that for $(x_1,....,x_n) \in E^n, \,\, n \geq k$ 

\begin{linenomath*}
\begin{equation*}
P(\textit{X}_1 = x_1,..,\textit{X}_n = x_n)= P( \textit{X}_1 = x_1,..,\textit{X}_k = x_k)
{\Pi}_{i=k + 1}^n Q(x_i \vert x_{i-k},..., x_{i-1})
\label{def_order}
\end{equation*}
\end{linenomath*}

for suitable transition probabilities $ Q(. \vert .)$. The class of processes that holds the above condition for a given
$k \geq 1$ will be denoted by ${\cal{M}}_k$, and ${\cal{M}}_0$ will denote the class of i.i.d. processes.
The \textbf{\textit{order of a Markov Chain}} in $\cup_{i=0}^\infty \,{\cal{M}}_i$ is the \textbf{\textit{smallest integer}} $\boxed{ \,\, \kappa \,\,}$ such that 
$\mathbb{X}= {\{ X_n \}}_{n \geq 0} \in {\cal{M}}_\kappa$.

\vspace{10pt}
Along the last few decades there has been a great number of research on the estimation of the order 
of a Markov Chains,
starting with  M.S. Bartlett \cite{Bartlett}, P.G. Hoel \cite{Hoel}, I.J. Good \cite{Good}, 
T.W. Anderson \& L.A. Goodman \cite{Anderson-Goodman}, P. Billingsley \cite{Billingsley1},
\cite{Billingsley2} among others, and more recently, H. Tong \cite{Tong}, 
G. Schwarz \cite{Schwarz}, R.W. Katz \cite{Katz}, I. Csiszar and P. Shields \cite{CsiszarShields_1}, 
L.C. Zhao et all \cite{Zhao} had contributed with new Markov chain order estimators.

Since 1973, H. Akaike \cite{Akaike}  entropic information criterion, known as AIC, has had a fundamental impact
in statistical model evaluation problems. The AIC has been applied by Tong, for example, to the 
problem of estimating the order of autoregressive processes, autoregressive integrated moving average 
processes, and Markov chains. The Akaike-Tong (AIC) estimator was  derived as an asymptotic approximate estimate of the 
Kullback-Leibler information discrepancy and provides a useful tool for evaluating models estimated by 
the maximum likelihood method. Later on,  Katz derived the asymptotic distribution of the 
estimator and showed  its inconsistency, proving that there is a positive probability of overestimating 
the true order no matter how large the sample size. 
Nevertheless, AIC is the most used and succesfull Markov chain  order estimator used at the present time, 
mainly because it is more  efficient than BIC  for small sample.

\vspace{30pt}
\section{Essentials on Some Estimators}
\subsection{Maximum Likelihood Methods}

The main consistent estimator alternative, the BIC estimator, does not perform too well for relatively 
small samples, 
as it was pointed out by Katz \cite{Katz} and Csiszar \& Shields \cite{CsiszarShields_1}. It is natural  to 
admit that the expansion of the Markov Chain complexity (size of the state space and order)  has significant 
influence on the sample size required for the identification of the unknown order, even though, most 
of the time it is difficult  to obtain  sufficiently large  samples.

Katz(1981) \cite{Katz} obtained the asymptotic distribution of ${\widehat \kappa}_{AIC}$ and proved 
its inconsistency showing 
the existence of a positive probability to overestimate the order. See also Shibata(1976) \cite{Shibata}. 
On the contrary Schwarz (1978) \cite{Schwarz} and Zhao(2001) \cite{Zhao} proved strong consistency for the estimators 
${\widehat \kappa}_{BIC}$ and ${\widehat \kappa}_{EDC}$, respectively. 

Clearly, for a given $\eta$,  $AIC(\eta)$ \cite{Tong},
$BIC(\eta)$ \cite{Schwarz} and $EDC(\eta)$ \cite{Zhao,Dorea} contain
much of the information concerning the sample's relative dependency, nevertheless  numerical simulations as well as 
theoretical considerations anticipates a great deal of variability for small samples.

\vspace{20pt}
Let $\textit{\textbf{X}}_1^n = (\textit{X}_1,...,\textit{X}_n)$  be a sample from a multiple stationary Markov chain 
$\mathbb{X}= {\{ \textit{X}_n \}_{n \geq 1}}$ 
of \textit{unknown order} ${\kappa}.$  

Assume that  $\mathbb{X}$ take values on a finite state space $E  = \{ 1,2,...,m\}$ with  transition probabilities  given by

\begin{linenomath*}
\begin{equation}
 p({x_{\kappa+1}\vert x_1^\kappa}) = P(\textit{X}_{n+1}=x_{n+1} \vert \textit{X}_{n-\kappa+1}^n=x_1^\kappa) > 0
\label{chain_cond_1}
\end{equation}
\end{linenomath*}

where $x_1^\kappa = (x_1,...,x_\kappa) = x_1^j\,x_{j+1}^\kappa  \in E^\kappa$.

\vspace{10pt}
Define 
\begin{linenomath*}
\begin{equation}
N(x_1^l \vert \textit{\textbf{X}}_1^n) = \sum_{j=1}^{n-l+1} 1(\textit{X}_j=x_1,...,\textit{X}_{j+l-1}=x_l)
\label{number_of_strings}
\end{equation}
\end{linenomath*}

i.e. the number of ocurrences of $x_1^l$ in $\textit{\textbf{X}}_1^n$. If $l=0$ \, we take $N(\,.\,\vert \textit{\textbf{X}}_1^n) = n$. 
The sums are taken over positive terms
$N(x_1^{l+1} \vert \textit{\textbf{X}}_1^n) > 0$, or else, we convention $0/0$ or $0.\infty$ as $0$.

\vspace{10pt}\
\begin{definition}
\begin{linenomath*}\label{def_X_alpha_ij}
For $a_1^\eta = (a_1,...,a_\eta) \in E^\eta$ and $j \in E$, let $\mathbb{X}_{a_1^\eta}$ be the empirical random variables,
extracted from the Markov chain sample $\textit{\textbf{X}}_1^n = (\textit{X}_1,...,\textit{X}_n)$  
\end{linenomath*}

\begin{linenomath*}
\begin{equation*}
\mathbb{X}_{a_1^\eta}\,:\, \textit{\textbf{X}}_1^n \longrightarrow 
\big( \mathbb{X}_{a_1^\eta}(1),...,\mathbb{X}_{a_1^\eta}(j),...,\mathbb{X}_{a_1^\eta}(m) \big) 
\end{equation*}
\end{linenomath*}

\begin{linenomath*}
\begin{equation}\label{probability_X_alpha_j}
\boxed{\,\,\,
\mathbb{X}_{a_1^\eta}(j) = 
\left(
\frac{N(a_1^\eta\,j \, \vert \, \textit{\textbf{X}}_1^n)}{N(a_1^\eta\,\vert\,\textit{\textbf{X}}_1^n)}
\right),
\,\,\,\,\,\, 1 \leq j \leq m \,\,\,}
\end{equation}
\end{linenomath*}

and

\begin{linenomath*}
\begin{equation*}
\mathbb{X}_{a_1^\eta}(i,j)\,:\, \textit{\textbf{X}}_1^n \longrightarrow 
\big( \mathbb{X}_{a_1^\eta \, i}(1),...,\mathbb{X}_{a_1^\eta \, i}(j),...,\mathbb{X}_{a_1^\eta \, i}(m) \big) 
\end{equation*}
\end{linenomath*}

\begin{linenomath*}
\begin{equation}\label{probability_X_alpha_ij}
\boxed{\,\,\,
\mathbb{X}_{a_1^\eta}(i,j) = \left(
\frac{N(a_1^\eta i\,j \, \vert \,  \textit{\textbf{X}}_1^n)}{N(a_1^\eta i\,\vert\,\textit{\textbf{X}}_1^n)}
\right),
\,\,\,\,\,\, 1 \leq i,j \leq m.\quad {\scriptstyle \blacklozenge} \,\,\,}
\end{equation}
\end{linenomath*}
\end{definition}

\vspace{30pt}
Let us define for the \textit{order} the log likelihood function

\begin{linenomath*}
\begin{equation*}
\log {\hat L}(\eta) = \sum_{j}^m {N(a_1^{\eta} \,\vert\,\textit{\textbf{X}}_1^n)} \left ( \sum_{a_1^{\eta}}
\frac{{N(a_1^{\eta}j\,\vert\,\textit{\textbf{X}}_1^n)}}{{N(a_1^{\eta}\,\vert\,\textit{\textbf{X}}_1^n)}} 
\log { \mathbb{X}_{a_1^\eta}(j) } \right ) 
\end{equation*}
\end{linenomath*}

\begin{linenomath*}
\begin{equation*}
\log {\hat L}(\eta) = \sum_{a_1^{\eta}} {N(a_1^{\eta} \,\vert\,\textit{\textbf{X}}_1^n)} \left ( \sum_{j}^m
\frac{{N(a_1^{\eta}j\,\vert\,\textit{\textbf{X}}_1^n)}}{{N(a_1^{\eta}\,\vert\,\textit{\textbf{X}}_1^n)}} 
\log { \mathbb{X}_{a_1^\eta}(j) } \right )       
\end{equation*}
\end{linenomath*}

\begin{linenomath*}
\begin{equation}\label{log_likelihood_function_1}
\boxed{\,
\log {\hat L}(\eta) = \sum_{a_1^{\eta}} \left [
{N(a_1^{\eta} \,\vert\,\mathbb{X}_1^n)} 
\left ( 
\sum_{j}^m  \mathbb{X}_{a_1^\eta}(j) \,\log { \mathbb{X}_{a_1^\eta}(j) }  
\right ) \right ]. 
\,}
\end{equation}
\end{linenomath*}

\vspace{10pt}
The estimators based on likelihood estimators and penalty functions, for Markov chains of order $\boxed{ \,\kappa \,}$ are defined, 
under the following hypothesis:
\begin{linenomath*}
\begin{equation*}
\boxed{ 
\textit{There exist a known}\,\, B \,\, \textit{so that} \,\, 0 \, \leq \kappa \, \leq B}
\end{equation*}
\end{linenomath*}
 
 as
\begin{linenomath*}
\begin{eqnarray}
 &&{\widehat \kappa}_{AIC} = \text{argmin}\{ AIC(\eta)\,;\, \eta=0,1,...,B\} \label{aic}\\
 &&{\widehat \kappa}_{BIC} = \text{argmin}\{ BIC(\eta)\,;\, \eta=0,1,...,B\} \label{bic}\\
 &&{\widehat \kappa}_{EDC} = \text{argmin}\{ EDC(\eta)\,;\, \eta=0,1,...,B\} \label{edc}
\end{eqnarray}
\end{linenomath*}

where

\begin{linenomath*}
\begin{eqnarray*}
  AIC(\eta) = -2 \log {\hat L}(\eta) \hspace{-20pt}&+& \hspace{-20pt} {\vert E \vert}^{\eta+1} \,\,2\,(\vert E \vert - 1), \\  
  BIC(\eta) = -2 \log {\hat L}(\eta) \hspace{-20pt}&+& \hspace{-20pt} {\vert E \vert}^{\eta+1} \,\,2\,(\vert E \vert - 1) 
                  \left( \frac{\log(n)}{2} \right ) \\
  EDC(\eta) = -2 \log {\hat L}(\eta) \hspace{-20pt}&+& \hspace{-20pt} {\vert E \vert}^{\eta+1} \,\,2\,(\vert E \vert - 1) 
                  \left(\frac{\log\log(n)}{2 (\vert E \vert - 1)}\right)\\
  AIC(\eta)  \,\, \leq   &EDC(\eta)&  \leq \,\, BIC(\eta).
\end{eqnarray*}
\end{linenomath*}

\vspace{35pt}
Finally, let us $\boxed{ \,\textit{\textbf{fix}} \,\,\, a_1^{\eta} \,}$ and consider the  function
$$\mathbb{L}_{[(n \,,\, a_1^{\eta})]} \,:\,{(0, 1)}^{m^2} \,\,\, \rightarrow \,\,\, \mathbb{R}^{+} $$
defined as:
\begin{linenomath*}
\begin{eqnarray}\label{function_aic_for_taylor}
&& \mathbb{L}_{[(n \,,\, a_1^{\eta})]}\Big(..., x(i,j), ...\Big) = {N(a_1^{\eta} \,\vert\, \textit{\textbf{X}}_1^n)} 
\Big( \sum_{i=1}^m  \mathbb{X}_{a_1^\eta}(i) \, \log \mathbb{X}_{a_1^\eta}(i) - \nonumber \\
&&\hspace{80pt} - 
\sum_{i=1}^m  \left[\frac{{N(a_1^{\eta} i \,\vert\, \textit{\textbf{X}}_1^n)}}{{N(a_1^{\eta} \,\vert\,\textit{\textbf{X}}_1^n)}}\right]
\sum_{j=1}^m  \mathbb{X}_{a_1^\eta}(i,j) \,\, \log \mathbb{X}_{a_1^\eta}(i,j) \Big). 
\end{eqnarray}
\end{linenomath*}

\vspace{25pt}
Later on in \textit{\textbf{Section 3.1}}, we shall analyse the behavior 
and  derivatives of $\textit{\textbf{L}}_{[(n \,,\, a_1^{\kappa})\,(i,k)]}$
which is just a generic representation of $\mathbb{L}_{[(n \,,\, a_1^{\kappa})\,(i,k)]}.$

$$\boxed{\,\, \textit{\textbf{L}}_{[(n \,,\, a_1^{\kappa})\,(i,k)]} \,\,:\,\, {(0,1)}  \rightarrow {\mathbb{R}}^{+} \,\,}$$

\begin{linenomath*}
\begin{equation} \label{log_likelihood_difference}
\boxed{\,\,\,
\textit{\textbf{L}}_{[( n \,,\, a_1^{\kappa})\,(i,k)]} 
\Big ( x(i,k) \Big )
= {\cal L}_{1} \Big( {\cal L}_{2} -  \pi(i) \,\, x(i,j) \, \log x(i,j) \Big )
\,\,\,}
\end{equation}
\end{linenomath*}

such that 
\begin{linenomath*}
\begin{equation*} 
\boxed{\,\,\,
\mathbb{L}_{[(n \,,\, a_1^{\kappa})\,(i,k)]}\Big(\mathbb{X}_{a_1^\kappa}(i,k)\Big) \equiv 
\textit{\textbf{L}}_{[(n \,,\, a_1^{\kappa})\,(i,k)]} \Big( x(i,j) \Big)
\,\,\,} 
\end{equation*}
\end{linenomath*}

where $\,{\cal L}_{1} = {N(a_1^{\kappa} \,\vert\, \textit{\textbf{X}}_1^n)}\,$ and 
$\,{\cal L}_{2} = \sum_{i=1}^m  \mathbb{X}_{a_1^\kappa}(i) \, \log \mathbb{X}_{a_1^\kappa}(i) \,$ are assumed \textit{\textbf{constants}} with 
respect to the  
the variables $x(i,j)$, with $x(i,j) = \mathbb{X}_{a_1^\kappa}(i,j)$ as in (\ref{probability_X_alpha_ij}), $\boxed{ \,\, \kappa \,\, }$ the Markov chain order
and
$$\pi(i) = \left [\frac{{N(a_1^{\kappa} i \,\vert\,\textit{\textbf{X}}_1^n)}}{{N(a_1^{\kappa} \,\vert\,\textit{\textbf{X}}_1^n)}}\right ],
\,\,\,\,\, 1 \leq i,j \leq m.$$ 

\vspace{30pt}
\subsection{$\chi^2$-divergence estimator}
We now briefly recall this new Markov  chain order's estimator referring the 
reader to (\cite{Baigorri}) for related details. 

\vspace{10pt}
\begin{definition}\label{global_dependency_level} 
Let $\textit{\textbf{X}}_n=\{X_i\}_{i=1}^n$ be a sample of a Markov chain $\mathbb{X}$ of order $\kappa \geq 0$,  $\mathbb{X}_{a_1^\eta}(i,j)$
as in (\ref{def_X_alpha_ij}), $\,\,\eta \geq 0 \,\,$ and ${{\varDelta}}_2( \mathbb{X}_{a_1^\eta}(i,j))$  the random 
variable defined as follows

\begin{linenomath*}
\begin{eqnarray*} 
&&{{\varDelta}}_2(\mathbb{X}_{a_1^\eta}(i,j)) = \\
&& = {N(a_1^{\eta} \,\vert\,\mathbb{X}_1^n)} \,\sum_{i=1}^m \sum\limits_{j=1}^m \left(  
\frac{ \Big( \mathbb{X}_{a_1^\eta}(i,j) - 
       \big[ \sum_{t=1}^m  \mathbb{X}_{a_1^\eta}(i,t) \big] \big[ \sum_{s=1}^m  \mathbb{X}_{a_1^\eta}(s,j) \big] \Big)^2}
     {( \sum_{t=1}^m \mathbb{X}_{a_1^\eta}(i,t) ) ( \sum_{s=1}^m  \mathbb{X}_{a_1^\eta}(s,j) )} \right) = \\
&& = {\cal G}_{(n,a_1^{\eta})} \,\sum_{i=1}^m \sum_{j=1}^m \mathbb{G}_{[(n,a_1^\eta)(i,j)]}.
\end{eqnarray*}
\end{linenomath*}
 
Assume that $V$ is  a $\chi^2$ random variable with $(m-1)^2$ degrees of freedom
where  ${\cal{P}}$ is the continuous strictly decreasing function ${\cal{P}}:{\mathbb{R}}^{+} \longrightarrow [0,1]$ 

\begin{linenomath*}
\begin{equation*}
{\cal{P}}(x) = P( V \geq x ), \,\,\, x \in {\mathbb{R}}^{+}.
\end{equation*}
\end{linenomath*}

The Local Dependency Level $\,\,\,{LDL}_n(a_1^\eta)\,\,\,$ and 
the Global Dependency Level $\,\,\,{GDL}_n(\eta)\,\,\,$, respectively, are defined as follows

\begin{linenomath*}
\begin{equation*}
{LDL}_n(a_1^\eta) =   
\frac{{{\varDelta}}_2( \mathbb{X}_{a_1^\eta}(i,j)}{2 \,\log(\log(n))},
\end{equation*}
\end{linenomath*}

\vspace{10pt}
\begin{linenomath*}
\begin{equation*}
{GDL}_n(\eta) = 
{\cal{P}} \left (\sum_{a_1^\eta \, \in E^\eta}\,\Big (\frac{N(a_1^\eta\,\vert\,\mathbb{X}_1^n)}{n}\Big )\,{LDL}_n(a_1^\eta) \right ).
\quad {\scriptstyle \blacklozenge}
\end{equation*}
\end{linenomath*}
\end{definition}

\vspace{10pt}
Finally, let us define the Markov chain order estimator based on the information contained in the vector $GDL_n.$

\begin{definition}
Given a fixed number $0 < B \in \mathbb{N}$, let us define the set ${\cal S} = \{0,1\}^{B+1}$ and the application 
$\,\, T\,:\,{\cal S} \rightarrow \mathbb{N}$

\begin{linenomath*}
\begin{eqnarray*}
T(s) &=& -1 \,\,\, \Leftrightarrow \,\,\, s_i=1, \,\,\, i=0,1,..,B \\
T(s) &=& \max_{0 \leq i \leq B} \left \{i \, : \, s_i=0, \, s_{i+1} = {\cal P(L)} \right \}, \,\,\, s=(s_0,s_1,...,s_B).
\quad {\scriptstyle \blacklozenge}
\end{eqnarray*}
\end{linenomath*}
\end{definition}

\vspace{10pt}
\begin{definition} 
Let $\textit{\textbf{X}}_1^n=\{X_i\}_{i=1}^n$ be 
a sample  for the Markov chain $\mathbb{X}$ of order $\kappa$, $0 \leq \kappa < B \in \mathbb{N}$  and  $\{{GDL}_n(i) \}_{i=1}^B$ as above. 
We define the order's estimator ${\kappa}_{GDL}(\mathbb{X}_1^n)$ as 

\begin{linenomath*}
\begin{equation*}
\widehat{\kappa}_{GDL}(\mathbb{X}_1^n) = T( \sigma_n ) + 1
\end{equation*}
\end{linenomath*}

with $\sigma_n \in {\cal S}$ so that $\,\, \forall \, s \in {\cal S}$

\begin{linenomath*}
\begin{equation*}
\sum_{i=0}^B \big ({GDL}_n(i) - \sigma_n(i) \big)^2 \leq  \sum_{i=0}^B ({GDL}_n(i) - s(i))^2. 
\quad {\scriptstyle \blacklozenge}
\end{equation*}
\end{linenomath*}
\end{definition}

\vspace{15pt}
Observe that  the \textit{\textbf{Local Dependency Level}} 
${LDL}_n(a_1^\eta)$  entirely relies on the just defined $\chi^2$-square  divergence estimator which itself
is the summation of several univariate random variables $\mathbb{G}_{[({n,\,a_1^\eta})\,(i,j)]}, \,\,\, 1 \leq i,j \leq m$

\begin{linenomath*}
\begin{eqnarray} \label{Delta_2_alpha_estimator}
&& \boxed{\,\,  \mathbb{G}_{[({n,\,a_1^\eta})\,(i,j)]} \,\,} =  \\
&& = \boxed{\,\,\sum_{i=1}^m \sum\limits_{j=1}^m   
\frac{ \Big( \mathbb{X}_{a_1^\eta}(i,j) - 
       \big[ \sum_{t=1}^m  \mathbb{X}_{a_1^\eta}(i,t) \big] \big[ \sum_{s=1}^m  \mathbb{X}_{a_1^\eta}(s,j) \big] \Big)^2}
     {\Big( \sum_{t=1}^m  \mathbb{X}_{a_1^\eta}(i,t) \Big) \Big( \sum_{s=1}^m  \mathbb{X}_{a_1^\eta}(s,j) \Big)}. \,\,} \nonumber                                            
\end{eqnarray}
\end{linenomath*}

\vspace{25pt}
Later on in \textit{\textbf{Section 3.2}}, we shall analyse the behavior of the deterministic function 
$\,\,\textit{\textbf{G}}_{[(n \,,\, a_1^{\kappa})\,(i,j)]}, \,\,\,\, 1 \leq i,j \leq m\,\,$ and their derivatives 

$$\textit{\textbf{G}}_{[(n \,,\, a_1^{\kappa})(i,j)]} \,:\, {(0,1)}^{2 m}  \rightarrow {\mathbb{R}}^{+}, \,\,\,\,\,\, 1 \leq i,j \leq m$$

with
\begin{linenomath*}
\begin{equation}\label{function_gdl_for_taylor}
\hspace{-8pt}\boxed{
\textit{\textbf{G}}_{[(n \,,\, a_1^{\kappa})\,(i,j)]}\big(...,x(s,t),...\big) = 
\frac{ \Big( x(i,j) - 
       \big[ \sum_{t=1}^m  x(i,t) \big] \big[ \sum_{s=1}^m  x(s,j) \big] \Big)^2}
     {\Big( \sum_{t=1}^m  x(i,t) ) ( \sum_{s=1}^m  x(s,j) \Big)} }
\end{equation}
\end{linenomath*}

such that 
\begin{linenomath*}
\begin{equation*} 
\boxed{\,\,\,
\mathbb{G}_{[(n \,,\, a_1^{\kappa})\,(i,j)]}\Big( ..., \mathbb{X}_{a_1^\kappa}(s,t),... \Big)  \equiv 
\textit{\textbf{G}}_{[(n \,,\, a_1^{\kappa})\,(i,j)]} \Big( ..., x(s,t),... \Big)
\,\,\,} 
\end{equation*}
\end{linenomath*}

with $x(s,t) = \mathbb{X}_{a_1^\kappa}(i,j), \,\,\, 1 \leq i,j \leq m$, \,\,  $\mathbb{X}_{a_1^\kappa}(i,j)$ as in (\ref{probability_X_alpha_ij}). 

\vspace{40pt}
\section{Deterministic Accessory Functions}\label{relevant_functions}
\subsection{Functions Related with AIC-Estimator}\label{log_likelihood_function}

Let us calculate the derivatives of the deterministic function $\textit{\textbf{L}}_{(n \,,\,a_1^{\kappa})}$ 
as in(\ref{log_likelihood_difference}), which for the sake of notational simplicity and for a fixed $n$ and $a_1^{\kappa}$, 
we'll \textit{\textbf{temporarily rename}} the function 

$$\boxed{\,\, \textit{\textbf{L}} = \textit{\textbf{L}}_{[(n \,,\, a_1^{\kappa})\,(i,k)]}, \,\,\,\,\,
\textit{\textbf{L}}\,:\, D_{\textit{\textbf{L}}} \subseteq (0,1)  \rightarrow \mathbb{R}, \quad 
\,\,}$$ 

$$ \boxed{\,\,\, D_\textit{\textbf{L}} = 
\{x(i,k) \,:\, x(i,k) \in (0,1)\},\,\,\,}$$ 

\begin{linenomath*} 
\begin{equation} \label{atomic_component_function_log_likelihood}
\boxed{\boxed{\textit{\textbf{L}} (\textit{\textbf{x}}) =  x{(i,k)} \log(x{(i,k)}). \,\,\,}}
\end{equation}
\end{linenomath*}

\vspace{30pt}
$\textbf{\textit{First Order Derivatives\,:}}$  

\begin{linenomath*} 
\begin{equation*}
\boxed{ \,\,\,
\frac{\partial \textit{\textbf{L}}}{\partial x{(i,k)}}  = 1 + \log(x{(i,k)}), \,\,\, 1 \leq i,k \leq m-1. \,\,\,}
\end{equation*}
\end{linenomath*}

\vspace{30pt}
$\textbf{\textit{Second Order Derivatives\,:}}$   
\begin{linenomath*} 
\begin{equation}\label{caso_2_L_parcial_x(i,k)_parcial_x(i,k)}
\boxed{ \,\,\,
\frac{\partial^2 \textit{\textbf{L}}}{\partial x^2{(i,k)} } = \frac{1}{x{(i,k)}},\,\,\, 1 \leq i,k \leq m-1, \,\,\,}
\end{equation}
\end{linenomath*}
\begin{linenomath*} 
\begin{equation}\label{caso_2_L_parcial_x(j,l)_parcial_x(i,k)}
\boxed{\,\,\, 
\frac{\partial^2 \textit{\textbf{L}}}{\partial x{(j,l)} \partial x{(i,k)}} = 0,\,\,\, 1 \leq i,k,l \leq m-1, \,\,\,} 
\end{equation}
\end{linenomath*}
respectively.

\vspace{30pt}
Later on we shall obtain the gradient vector  and the hessian matrix
$$\nabla_{{\textit{\textbf{L}}}}\left(\varLambda_{a_1^\kappa}(o)\right)
\, , \, {\cal H}_{\textit{\textbf{L}}}\left(\varLambda_{a_1^\kappa}(o)\right)$$
of the function ${\textit{\textbf{L}}}$ at a point 
$${\varLambda}_{a_1^\kappa}(o)= \Big(..., E\big(\mathbb{X}_{a_1^\kappa}(i,k)\big),...\Big).$$

\vspace{50pt}
\subsection{Functions Related with GDL-Estimator}\label{local_dependency_functions}
Herein, we shall consider the deterministic multivariate set of functions, as in (\ref{function_gdl_for_taylor}),
\begin{linenomath*}
\begin{equation*} 
\boxed{\boxed{\,\,\,
\textit{\textbf{G}}_{[(n \,,\, a_1^{\kappa})\,(i,k)]}\Big(x{(i,k)},h(i),v(k) \Big) = \frac{\Big(x{(i,k)} - h(i) v(k) \Big)^2}{h(i) v(k)}, 
\,\,\, 1 \leq i,k \leq m.\,\,\,}}
\end{equation*}
\end{linenomath*}
which, after fixing $(i,k)$, we \textbf{\textit{temporarily rename}} it  as follows:  

\vspace{10pt}
$$\boxed{\,\,\,  \textit{\textbf{G}} \equiv \textit{\textbf{G}}_{[(n \,,\, a_1^{\kappa})\,(i,k)]}, \,\,\,\,\,\,\,\,\,\,
           \textit{\textbf{G }} : D_{\textit{\textbf{G}}} \subseteq (0,1)^3  \rightarrow \mathbb{R}, \,\,\, }$$

\begin{linenomath*}
\begin{equation*} 
\boxed{\,\,\, D_{\textbf{G}} = \{  \textit{\textbf{x}} \in (0,1)^3  \,:\,  \textit{\textbf{x}}= ( x, h, v )\}, \,\,\,}  
\end{equation*}
\end{linenomath*}

\begin{linenomath*} 
\begin{eqnarray} \label{atomic_component_function_G_d_l}
&& \boxed{\,\,\, \textit{\textbf{G}} (\textit{\textbf{x}})\, = \,\frac{( x -  h v  )^2}{ h v }.\,\,\, }
\end{eqnarray}
\end{linenomath*}

\vspace{25pt}
\textbf{\textit{First Order Derivatives\,\,:\,}} 

\begin{linenomath*} 
\begin{equation*}\label{G_first_derivatives}
\boxed{ \,\,\, \frac{\partial \textit{\textbf{G}}}{\partial x}  =  
\frac{2 x}{h v} - 2, \,\,\,\,\,\,\, \frac{\partial \textit{\textbf{G}}}{\partial h} = \frac{- x^2}{v h^2} + v, \,\,\,\,\,\,\,
\frac{\partial \textit{\textbf{G}}}{\partial v} = \frac{- x^2}{h v^2} + h.\,\,\,}
\end{equation*}
\end{linenomath*}

\vspace{35pt}
\textbf{\textit{Second Order Derivatives\,\,:\,}} 
           
\begin{linenomath*} 
\begin{equation}\label{G_second_derivatives}
\boxed{\frac{\partial^2 \textit{\textbf{G}}}{\partial x^2} = \frac{2}{h v}, \,\,\,\,\,\,\,
\,\,   \frac{\partial^2 \textit{\textbf{G}}}{\partial h^2} = \frac{- x^2}{h^3  v},  \,\,\,\,\,\,\,
\,\,   \frac{\partial^2 \textit{\textbf{G}}}{\partial v^2} = \frac{- x^2}{v^3 h},
 \,\,}  
\end{equation}
\end{linenomath*}

\begin{linenomath*} 
\begin{equation}\label{G_second_derivatives}
\boxed{\,\,
\,\,\frac{\partial^2 \textit{\textbf{G}}}{\partial h \partial x} = \frac{- 2 x}{h^2 v},\,\,\,\,\,\,\,  
\,\,\frac{\partial^2 \textit{\textbf{G}}}{\partial v \partial x} = \frac{- 2 x}{h v^2}, \,\,\,\,\,\,\,
\,\,\frac{\partial^2 \textit{\textbf{G}}}{\partial v \partial h} = \frac{ x^2}{h^2 v^2} + 1. 
\,\,}
\end{equation}
\end{linenomath*}

\vspace{30pt}
Likewise, as in the previous subsection we get the gradient vector  and the hessian matrix
$\,\, \nabla_{{\textit{\textbf{G}}}}\,\left( \Gamma_{a_1^\kappa}(o) \right) \, , 
   \, {\cal H}_{\textit{\textbf{G}}}\,\left( \Gamma_{a_1^\kappa}(o) \right) \,\,$
of the function ${\textit{\textbf{G}}}(x,h,v)$ at a point 
$$\Gamma_{a_1^\kappa}(o)= \Big(\,E(\mathbb{X}_{a_1^\kappa}(i,k)), E(\mathbb{H}_{a_1^\kappa}(i)), E(\mathbb{V}_{a_1^\kappa}(k))\,\Big).$$

\vspace{30pt}
\section{Multivariate Variances}

Focusing on  $\textit{\textbf{G}}_{[(n \,,\, a_1^{\kappa})\,(i,k)]} \,for \,\, fixed \,\, (i,k)$,  
$\kappa$ the order of the Markov chain and 
$\textit{\textbf{G}}_{[(n \,,\, a_1^{\kappa})\,(i,k)]}$ as in (\ref{function_gdl_for_taylor}),
let us recall the empirical random variables, introduced in Definition \ref{def_X_alpha_ij}

$$ \mathbb{X}_{a_1^\kappa}{(i,k)} = 
\frac{N(i\,a_1^\kappa\,k \, \vert \, \textit{\textbf{X}}_1^n)}{N(a_1^\kappa\,\vert\, \textit{\textbf{X}}_1^n))}, \,\,\,\, 
\mathbb{H}_{a_1^\kappa}{(i)} = \sum_{t=1}^m \mathbb{X}_{a_1^\kappa}{(i,t)}, 
\,\,\,\, \mathbb{V}_{a_1^\kappa}{(k)} = \sum_{s=1}^m \mathbb{X}_{a_1^\kappa}{(s,k)}.$$

Observe that the Markov chain we are interested in, has order $\kappa$ and it is clear that 
\begin{linenomath*} 
\begin{equation*}
\boxed{\,\,\,
\mathbb{X}_{a_1^\kappa}{(i,k)} \,\,\, \text{\footnotesize independent random variables}, \,\,\,\, 1 \leq i,k \leq m, \,\,\,}
\end{equation*}
\end{linenomath*}

\begin{linenomath*} 
\begin{equation}\label{G_(i,k)_independence}
\boxed{\,\, \textit{\textbf{G}}_{[(n \,,\, a_1^{\kappa})\,(i,j)]}, \,\,\,\,\, 1 \leq i,k \leq m-1, \,\,}
\end{equation}
\end{linenomath*}

are  independent, with
\begin{linenomath*} 
\begin{eqnarray*}
\textit{\textbf{G}}_{[(n \,,\, a_1^{\kappa})\,(i,j)]} 
\left(\mathbb{X}_{a_1^\kappa}{(i,k)}, \mathbb{H}_{a_1^\kappa}{(i)}, \mathbb{V}_{a_1^\kappa}{(k)} \right)
= \frac{\left(\mathbb{X}_{a_1^\kappa}{(i,k)} -\mathbb{H}_{a_1^\kappa}{(i)}\,\mathbb{V}_{a_1^\kappa}{(k)}\right)^2}
         {\mathbb{H}_{a_1^\kappa}{(i)} \mathbb{V}_{a_1^\kappa}{(k)}}
\end{eqnarray*}
\end{linenomath*}
as well as for \textbf{\textit{adequatly}} sample size $\textit{\textbf{n}}$  the random variables 

\begin{linenomath*} 
\begin{equation*}
\boxed{\,\,\,
E\big ( \mathbb{H}_{a_1^\kappa}{(i)}  \big ) \approxeq \frac{1}{m},  
\,\,\,\,\,\, E\big ( \mathbb{V}_{a_1^\kappa}{(k)}  \big ) \approxeq m\,E\big ( \mathbb{X}_{a_1^\kappa}{(s,k)}  \big ), \,\,\,\, 1 \leq s \leq m. \,\,\,}
\end{equation*}
\end{linenomath*}

\vspace{20pt}
For the sake of notation's simplicity, we temporarily rename 
$$\textit{\textbf{G}}_{[(n \,,\, a_1^{\kappa})\,(i,j)]} \big( x{(i,k)}, h(i), v(k) \big)$$    
as $\textit{\textbf{G}}(x, h, v)$  where its derivatives, as well as the variances, covariances and related  information of 
$\{ \, \mathbb{X}_{a_1^\kappa}{(i,k)},\,\mathbb{H}_{a_1^\kappa}{(i)}, \,\mathbb{V}_{a_1^\kappa}{(k)}\, \}$ shall be as follows:

\begin{linenomath*} 
\begin{equation*}
\boxed{\,\,\,
E\left( \frac{\mathbb{X}_{a_1^\kappa}{(i,k)}}{\mathbb{H}_{a_1^\kappa}(i) \mathbb{V}_{a_1^\kappa}(k)} \right) 
\,\approxeq\, \frac{E\left( \mathbb{X}_{a_1^\kappa}{(i,k)}\right)} {E\left( \mathbb{H}_{a_1^\kappa}(i)\right) E\left( \mathbb{V}_{a_1^\kappa}(k)\right)} 
\,\approxeq\, \frac{ E(\mathbb{X}_{a_1^\kappa}{(i,k)})}{ \frac{1}{m} \big( m \,E(\mathbb{X}_{a_1^\kappa}{(i,k)} )\big)  } \,\approxeq\, 1,
\,\,\, }
\end{equation*}
\end{linenomath*}

%

\vspace{35pt}
\begin{linenomath*} 
\begin{eqnarray*}
&&\hspace{-40pt}\boxed{ \,\,
\sum_{\alpha,\beta \in \{x,h,v\}} \frac{\partial \textit{\textbf{G}}}{\partial \alpha}(x,h,v) 
         \frac{\partial \textit{\textbf{G}}}{\partial \beta }(x,h,v) \,\, } \approxeq\,\, \\
&\approxeq& \left( 1 - \Big(\frac{x}{h\,v} \Big)^2 \right) \Big( 4 + h^2 + v^2 \Big) 
+ \frac{2 x}{h}\left( \Big(\frac{x}{h\,v}\Big) - \Big(\frac{x}{h\,v}\Big)^2 \right) + \\
&&\frac{2 x}{v}\left( \Big(\frac{x}{h\,v}\Big) - \Big(\frac{x}{h\,v}\Big)^2 \right) + 
  \,x \left( \Big(\frac{x}{h\,v}\Big)^3 - 2 \,\Big(\frac{x}{h\,v}\Big)^2 + \Big(\frac{x}{h v}\Big)  \right) \\
&\approxeq& \boxed{\,\,0\,\,.}
\end{eqnarray*}
\end{linenomath*}
and
\begin{linenomath*} 
\begin{equation*}
\hspace{-30pt}\boxed{
\textbf{cov}_{(\mathbb{X}_{a_1^\kappa}{(i,k)}\,,\,\mathbb{H}_{a_1^\kappa}(i))} \,\approxeq\, 
\textbf{cov}_{(\mathbb{X}_{a_1^\kappa}{(i,k)}\,,\,\mathbb{V}_{a_1^\kappa}(k))} \,\approxeq\,
\textbf{cov}_{(\mathbb{H}_{a_1^\kappa}(i)\,,\,\mathbb{V}_{a_1^\kappa}(k))} \,\approxeq\,
\sigma^2_{\mathbb{X}_{a_1^\kappa}{(i,k)}}.
}
\end{equation*}
\end{linenomath*}

\vspace{25pt}
Likewise,
\begin{linenomath*} 
\begin{equation*}
\sum_{\alpha, \beta, \gamma \in \{x,h,v\}} 
\frac{\partial \textit{\textbf{G}}}{\partial \alpha}(x,h,v) \frac{\partial^2 \textit{\textbf{G}}}{\partial \beta \partial \gamma }(x,h,v) 
\,\,\approxeq\,\, 0.
\end{equation*}
\end{linenomath*}

\begin{linenomath*} 
\begin{equation*}
\boxed{\,\,\,
\textbf{cov}_{(\mathbb{X}_{a_1^\kappa}{(i,k)}\,,\,\mathbb{X}_{a_1^\kappa}^2{(i,k)})} \,\approxeq\, 
\textbf{cov}_{(\mathbb{X}_{a_1^\kappa}{(i,k)}\,,\,\mathbb{H}_{a_1^\kappa}^2(i))} \,\approxeq\, 
\textbf{cov}_{(\mathbb{X}_{a_1^\kappa}{(i,k)}\,,\,\mathbb{V}_{a_1^\kappa}^2(k))} \,\approxeq\,
\sigma^3_{{\mathbb{X}_{a_1^\kappa}{(i,k)}}}
\,\,\,}
\end{equation*}
\end{linenomath*}

\begin{linenomath*} 
\begin{equation*}
\boxed{\,\,\,
\textbf{cov}_{(\mathbb{X}_{a_1^\kappa}{(i,k)}\,,\,\mathbb{H}_{a_1^\kappa}^2(i))} \,\approxeq\, 
\textbf{cov}_{(\mathbb{X}_{a_1^\kappa}{(i,k)}\,,\,\mathbb{V}_{a_1^\kappa}^2(k))} \,\approxeq\, 
\sigma^3_{{\mathbb{X}_{a_1^\kappa}{(i,k)}}}. \,\,\,}
\end{equation*}
\end{linenomath*}

\vspace{30pt}
Finally,
\begin{linenomath*} 
\begin{eqnarray*}
&& \boxed{\,\, \sum_{\alpha,\beta \in \{x,h,v\}} 
               \frac{\partial^2 \textit{\textbf{G}}}{\partial \alpha^2}(x,h,v) 
               \frac{\partial^2 \textit{\textbf{G}}}{\partial  \beta^2}(x,h,v) \,}
                                              \,\,\approxeq\,\, 
             \left(\frac{\partial^2 \textit{\textbf{G}}}{\partial x^2}\right)^2 +
                   \frac{\partial^2 \textit{\textbf{G}}}{\partial x^2} \frac{\partial^2 \textit{\textbf{G}}}{\partial v^2} +
                   \frac{\partial^2 \textit{\textbf{G}}}{\partial x^2} \frac{\partial^2 \textit{\textbf{G}}}{\partial v \partial x}  + \\
      && +   \left(\frac{\partial^2 \textit{\textbf{G}}}{\partial v^2}\right)^2  + 
                   \frac{\partial^2 \textit{\textbf{G}}}{\partial v^2} \frac{\partial^2 \textit{\textbf{G}}}{\partial v \partial x} +
             \left(\frac{\partial^2 \textit{\textbf{G}}}{\partial v \partial x}\right)^2 = 
                   \frac{1}{h^2 v^2} \Big( 4 + \frac{x^4}{v^4} + \frac{4 x^3}{v^3} - \frac{8 x}{v} \Big)  \,\,\approxeq\,\, \\
      && \,\,\approxeq\,\, \frac{1}{x^2} \Big( 4 + \frac{1}{m^4} + \frac{4 }{m^3} - \frac{8} {m} \Big) \,\,\approxeq\,\,
         \frac{4}{x^2}\big( 1 - \frac{2}{m} \big) + \frac{1}{m^3} + \frac{1}{m^4}\,\,\approxeq\,\,\boxed{\,\,\frac{4}{x^2}\,\,}
\end{eqnarray*}
\end{linenomath*}

\begin{linenomath*} 
\begin{equation*}
\boxed{\,\,\,
\textbf{cov}_{(\mathbb{X}_{a_1^\kappa}^2{(i,k)}\,,\,\mathbb{X}_{a_1^\kappa}^2{(i,k)})} \,\approxeq\,
\textbf{cov}_{(\mathbb{X}_{a_1^\kappa}^2{(i,k)}\,,\,\mathbb{H}_{a_1^\kappa}^2(i))} \,\approxeq\,
\textbf{cov}_{(\mathbb{X}_{a_1^\kappa}^2{(i,k)}\,,\,\mathbb{V}_{a_1^\kappa}^2(k))} \,\approxeq\, 
\sigma^4_{{\mathbb{X}_{a_1^\kappa}{(i,k)}}} \,\,\,}
\end{equation*}
\end{linenomath*}

\begin{linenomath*} 
\begin{equation*}
\boxed{\,\,\,
\textbf{cov}_{(\mathbb{H}_{a_1^\kappa}^2(i)\,,\,\mathbb{H}_{a_1^\kappa}^2(i))} \,\approxeq\,
\textbf{cov}_{(\mathbb{H}_{a_1^\kappa}^2(i)\,,\,\mathbb{V}_{a_1^\kappa}^2(k))} \,\approxeq\,
\textbf{cov}_{(\mathbb{V}_{a_1^\kappa}^2(k)\,,\,\mathbb{V}_{a_1^\kappa}^2(k))} \,\approxeq\, 
\sigma^4_{{\mathbb{X}_{a_1^\kappa}{(i,k)}}}.
\,\,\,}
\end{equation*}
\end{linenomath*}

\vspace{30pt}
Let us denote by  ${\cal B} \in {R}^3$, the unit ball  centered at the point

$$\varGamma_{a_1^\kappa}(o) = 
\Big(\,E(\mathbb{X}_{a_1^\kappa}{(i,k)}),E(\mathbb{H}_{a_1^\kappa}{(i)}), E(\mathbb{V}_{a_1^\kappa}{(k)}\,\Big)$$

with
$$\omega = \Big(\mathbb{X}_{a_1^\kappa}(i,k),\mathbb{H}_{a_1^\kappa}(i),\mathbb{V}_{a_1^\kappa}(k)\Big), 
\,\,\, \varDelta = \omega - \varGamma_{a_1^\kappa}(o).$$ 

\vspace{20pt}
Taylor (\cite{Apostol}) showed that there exist $ 0 < c_g, c_l < 1$ such that
\vspace{25pt}
\begin{linenomath*}
\begin{eqnarray*}
\textit{\textbf{G}}_{[(n \,,\, a_1^{\kappa})\,(i,j)]}(\omega) &=& \\
 &&\hspace{-40pt} \textit{\textbf{G}}_{(n \,,\, a_1^{\kappa})}[i,j]\Big(\varGamma_{a_1^\kappa}(o)\Big) + 
         \nabla_{\textit{\textbf{G}}_{(n \,,\, a_1^{\kappa})}[i,j]} \Big(\varGamma_{a_1^\kappa}(o)\Big)\,.(\omega - \varGamma_{a_1^\kappa}(o)) + \\
\hspace{-25pt}&& + \frac{1}{2!} \, (\omega - \varGamma_{a_1^\kappa}(o)) \,.\, {\cal H}_{\textit{\textbf{G}}_{(n \,,\, a_1^{\kappa})}[i,j]}
\Big( \varGamma_{a_1^\kappa}(o) + c_g \varDelta \Big) \,.\, (\omega - \varGamma_{a_1^\kappa}(o))^{\,t}         
\end{eqnarray*}
\end{linenomath*}

\vspace{15pt}
where the variance of 
$$\textit{\textbf{G}}_{[(n \,,\, a_1^{\kappa})\,(i,k)]} 
\Big( \mathbb{X}_{a_1^\kappa}{(i,k)},\,\mathbb{H}_{a_1^\kappa}{(i)},\,\mathbb{V}_{a_1^\kappa}{(k)} \Big)$$ 

is
\begin{linenomath*}
\begin{eqnarray*}
&&\sigma^2_{ \textit{\textbf{G}}_{[(n \,,\, a_1^{\kappa})\,(i,j)]}}= \\
&& = E \left[ \textit{\textbf{G}}_{(n \,,\, a_1^{\kappa})}[i,j]\big( \omega \big)                           
               - {\textit{\textbf{G}}_{(n \,,\, a_1^{\kappa})}[i,j]} \Big( \varGamma_{a_1^\kappa}(o) \Big) \right]^2 = \\                                    
&& = E \left[ \nabla_{\textit{\textbf{G}}_{(n \,,\, a_1^{\kappa})}[i,j]}\big(\varGamma_{a_1^\kappa}(o)\big)
 \,.\,\varDelta + \,\frac{1}{2!} \varDelta \,.\, {\cal H}_{\textit{\textbf{G}}_{(n \,,\, a_1^{\kappa})}[i,j]} 
               \Big( \varGamma_{a_1^\kappa}(o) + c_g \,(\omega - \varGamma_{a_1^\kappa}(o) ) \Big) \,.\, \varDelta^t \right ]^2.
\end{eqnarray*}
\end{linenomath*}

\vspace{10pt}
\begin{linenomath*}
\begin{equation*}\label{G(i,k)_variance}
\boxed{\,\,\,
\sigma^2_{\mathbb{G}{[(n,\, a_1^{\kappa})\,(i,k)]}}
  \,\approxeq\, 
\frac{1}{\big[2!]^2} \,\varPhi''_g \Big( \varGamma_{a_1^\kappa}(o) + c_g \,(\omega - \varGamma_{a_1^\kappa}(o) ) \Big) 
\, \sigma^4_{\mathbb{X}_{a_1^\kappa}{(i,k)}}\,\,\,}
\end{equation*}
\end{linenomath*}

\vspace{20pt}
where $\omega = \Big(\mathbb{X}_{a_1^\kappa}(i,k),\mathbb{H}_{a_1^\kappa}(i),\mathbb{V}_{a_1^\kappa}(k)\Big) \in {\cal B}$ and 

\begin{linenomath*}
\begin{eqnarray*}
\varPhi''_g \Big( x(i,k),h(i),v(k) \Big) &=& \left(\frac{1}{h^2(i) v^2(k)} \right) \left[ 4  \,+\, 
                    \frac{x^4(i,k)}{v^4(k)}\,+\, 4 \,\frac{x^3(i,k)}{v^3(k)}\,-\, 8 \,\frac{x(i,k)}{v(k)} \right]\\
                    &\approxeq& \, \frac{4}{x(i,k)}
\end{eqnarray*}
\end{linenomath*}

with variance for $ 1 \leq i,k \leq m-1 $
\begin{linenomath*}
\begin{equation*}
\boxed{\,\,\,
\sigma^2_{\mathbb{G}{[(n,\, a_1^{\kappa})\,(i,k)]}}
  \,\approxeq\, 
\frac{\sigma^4_{\mathbb{X}_{a_1^\kappa}{(i,k)}}}
     {\left[ \varGamma_{a_1^\kappa}(o) + c_g \,(\omega - \varGamma_{a_1^\kappa}(o) )  \right]^2}.
\,\,\,}
\end{equation*}
\end{linenomath*}

\vspace{20pt}
and, by (\ref{G_(i,k)_independence}), the total variance

\begin{linenomath*}
\begin{equation*}\label{G_total_variance}
\boxed{\,\,
\sigma^2_{ \mathbb{G}_{[({n,\,a_1^\eta})]} }= 
\sum_{i=1}^{m-1} \sum_{k=1}^{m-1} \sigma^2_{\mathbb{G}{[(n \,,\, a_1^{\kappa})\,(i,k)]}}. 
\,\,}
\end{equation*}
\end{linenomath*}

\vspace{40pt}
Exactly as before we can obtain the total variance of $\mathbb{L}_{(n, a_1^\kappa)}\big(\mathbb{X}_{a_1^\kappa}(i,k)\big)$ 
and defining 

\begin{linenomath*}
\begin{eqnarray*}
&&\sigma^2_{\mathbb{L}_{[(n,\,a_1^\kappa)\,(i,k)]}} = \\
&& = E \Big[ \mathbb{L}_{[(n, a_1^\kappa)\,(i,k)]} \Big( \mathbb{X}_{a_1^\kappa}(i,k) \Big)                           
           - \mathbb{L}_{[(n, a_1^\kappa)\,(i,k)]} \Big( E(\mathbb{X}_{a_1^\kappa}(i,k)) \Big) \Big ]^2 = \\                                
&& = E \Big [ \nabla_{\mathbb{L}_{[(n, a_1^\kappa)\,(i,k)]}}.\Big( \mathbb{X}_{a_1^\kappa}(i,k) - E(\mathbb{X}_{a_1^\kappa}(i,k)) \Big) + \\                                 
&& \hspace{20pt} + \,\frac{1}{2!} \, {\cal H}_{ \mathbb{L}_{[(n, a_1^\kappa)\,(i,k)]} } 
  \Big( E(\mathbb{X}_{a_1^\kappa}) + c_l \big[x(i,k) - E(\mathbb{X}_{a_1^\kappa})\big] \Big) \,.\, 
               \Big( \mathbb{X}_{a_1^\kappa}(i,k) - E(\mathbb{X}_{a_1^\kappa}(i,k)) \Big)^2 \Big ]^2.
\end{eqnarray*}
\end{linenomath*}

\begin{linenomath*}
\begin{equation*}\label{L(i,k)_variance}
\boxed{
\sigma^2_{{\mathbb{L}}_{[(n,\,a_1^\kappa)(i,k)]}}  \,\approxeq\,
\left[ 1 + \ln(x{(i,k)}) \right]^2  \, \sigma^2_{\mathbb{X}_{a_1^\kappa}{(i,k)}} +
 \, \varPhi''_l\big(x(i,k)\big)\,\, \sigma^4_{\mathbb{X}_{a_1^\kappa}{(i,k)}},\, 1 \leq i,k \leq m 
}  
\end{equation*} 
\end{linenomath*}

\vspace{20pt}
where $\omega = \big( x(i,k),h(i),v(k) \big) \in {\cal B}$,   

\begin{linenomath*}
\begin{equation*}
\varPhi''_l\big(x(i,k)\big) =  \left[ \frac{1}{x^2(i,k)} \right] \\
\end{equation*} 
\end{linenomath*}

with variance for $ 1 \leq i,k \leq m-1 $
\begin{linenomath*}
\begin{equation*}
\boxed{\,\,\,
\sigma^2_{{\mathbb{L}}_{[(n, a_1^\kappa)(i,k)]}}  \,\,\approxeq\,\, 
\left[ 1 + \ln(x{(i,k)}) \right]^2  \, \sigma^2_{\mathbb{X}_{a_1^\kappa}{(i,k)}} + 
\frac{\sigma^4_{\mathbb{X}_{a_1^\kappa}{(i,k)}}}{  \left[ \Big(E(\mathbb{X}_{a_1^\kappa}) + c_l \big[x(i,k) - E(\mathbb{X}_{a_1^\kappa})\big] \Big) \right]^2 },
\,\,\,}
\end{equation*}
\end{linenomath*}

\vspace{20pt}
and, by (\ref{G_(i,k)_independence}), the total variance

\begin{linenomath*}
\begin{equation*}\label{L_total_variance}
\boxed{\,\,\,\sigma^2_{{\mathbb{L}}_{[(n,\,a_1^\kappa)]}} \,\approxeq\,\, 
\sum_{k=1}^{m-1} \sum_{i=1}^{m-1} \sigma^2_{{\mathbb{L}}_{(n,\,a_1^\kappa)}{(i,k)}}.\,\,\,} 
\end{equation*}
\end{linenomath*}

\vspace{30pt}
\section{Conclusion}
The purpose of this work was the comparative analysis  of the non asymptotic behavior for the 
estimators {\textit{AIC}}(\ref{aic}), {\textit{BIC}}(\ref{bic}), {\textit{EDC}}(\ref{edc}), versus the  estimator defined in 
Definition \ref{global_dependency_level} and named as \textit{Global Depency Level-{GDL}}, for details see (\cite{Baigorri}).

The  \textit{GDL} uses a function different to the log likelihoog function applied to 
the  sample, which makes the estimator perform  in a quite different form. It is strongly consistent and  
more efficient than \textit{ AIC}(inconsistent),  outperforming the well established and consistent \textit{BIC} and \textit{EDC,} 
mainly on reasonable small samples.

The estimators  just mentioned are based on the composition of the empirical random variables with two diferent \textit{deterministic} 
functions. The log likelihood approach, as in (\ref{log_likelihood_difference}), with

\begin{linenomath*}
\begin{equation*}
\boxed{\,\,\,
\textit{\textbf{L}}_{[(n \,,\, a_1^{\kappa})\,(i,k)]} 
\Big ( x(i,k) \Big )
= {\cal L}_{1} \Big( {\cal L}_{2} -  \pi(i) \, x(i,k) \, \log x(i,k) \Big )
\,\,\,}
\end{equation*}
\end{linenomath*}

\vspace{10pt}

or, the GDL approach, as in (\ref{Delta_2_alpha_estimator}), with

\begin{linenomath*}
\begin{equation*}
\boxed{\,\,
\textit{\textbf{G}}_{[(n \,,\, a_1^{\kappa})\,(i,j)]}\Big ( ..., x(i,j),... \Big ) = 
\frac{ \Big( x(i,j) - 
       \big[ \sum_{t=1}^m  x(i,t) \big] \big[ \sum_{s=1}^m  x(s,j) \big] \Big)^2}
     {( \sum_{t=1}^m  x(i,t) ) ( \sum_{s=1}^m  x(s,j) )}. \,\,}
\end{equation*}
\end{linenomath*}

\vspace{30pt}
Since the sample only depends on the Markov chain $ \mathbb{X}_1^n $ and its size $\textbf{n}$, once the sample is chosen, 
the entirely responsibles for the estimator's variance are the following random variables:

\begin{linenomath*}
\begin{equation*}
\boxed{\,\,\,
\mathbb{L}_{[(n \,,\, a_1^{\kappa})(i,k)]} 
= {\cal L}_{1} \Big( {\cal L}_{2} - \pi(i) \,\, \mathbb{X}_{a_1^\kappa}(i,k) \log \mathbb{X}_{a_1^\kappa}(i,k)
\,\,\,}
\end{equation*}
\end{linenomath*}

and 
\begin{linenomath*}
\begin{equation*}
\hspace{-8pt}\boxed{\,\,\,
\mathbb{G}_{[(n \,,\, a_1^{\kappa})\,(i,k)]} = 
\sum_{k=1}^{m}\sum_{i=1}^{m}\frac{ \Big( \mathbb{X}_{a_1^\kappa}(i,k) - 
       \big[ \sum_{t=1}^m  \mathbb{X}_{a_1^\kappa}(i,t) \big] \big[ \sum_{s=1}^m  \mathbb{X}_{a_1^\kappa}(s,k) \big] \Big)^2}
     {( \sum_{t=1}^m  \mathbb{X}_{a_1^\kappa}(i,t) ) ( \sum_{s=1}^m  \mathbb{X}_{a_1^\kappa}(s,k) )} \,\,\,}
\end{equation*}
\end{linenomath*}

\vspace{20pt}
with variances for $ 1 \leq i,k \leq m-1 $
\begin{linenomath*}
\begin{equation*}
\boxed{\boxed{\,\,\,
\sigma^2_{{\mathbb{L}}_{[(n, a_1^\kappa)(i,k)]}}  \,\,\approxeq\,\, 
\left[ 1 + \ln(x{(i,k)}) \right]^2  \, \sigma^2_{\mathbb{X}_{a_1^\kappa}{(i,k)}} + 
\frac{\sigma^4_{\mathbb{X}_{a_1^\kappa}{(i,k)}}}{  \left[ \Big(E(\mathbb{X}_{a_1^\kappa}) + c_l \big[x(i,k) - E(\mathbb{X}_{a_1^\kappa})\big] \Big) \right]^2 },
\,\,\,}}
\end{equation*}
\end{linenomath*}

and
\begin{linenomath*}
\begin{equation*}
\boxed{\boxed{\,\,\,
\sigma^2_{\mathbb{G}{[(n,\, a_1^{\kappa})\,(i,k)]}}
  \,\approxeq\, 
\frac{4\,\,\sigma^4_{\mathbb{X}_{a_1^\kappa}{(i,k)}}}
     {\left[ \varGamma_{a_1^\kappa}(o) + c_g \,(\omega - \varGamma_{a_1^\kappa}(o) )  \right]^2}.
\,\,\,}}
\end{equation*}
\end{linenomath*}
respectively.

\vspace{15pt}
Finally the reader should notice that the log likelihood based estimators are heavily affected by $\boxed{\log( x{(i,k))} )}$ which
in  cases where the Markov chain intrisically presents  empirical random variables $\textit{\textbf{X}}_{a_1^\kappa}(i,k)$ with 
small expectations,
the fluctuating  values of $x(i,k)$ converging to $E(\textit{\textbf{X}}_{a_1^\kappa}{(i,k)}) \,\approxeq\, 0$ imposes the coefficients 
$\left[ 1 + \log(x{(i,k)}) \right]^2$ and its variance $\sigma^2_{{\textit{\textbf{L}}}_{(n, a_1^\kappa)}}$ a great deal of instability or variance.

The following Appendix presents a few examples exhibiting such anomaly.

\vspace{30pt}
\section{Appendix}
\subsection{Numerical Evidence}
In what follows we shall compare the non-asymptotic performance, mainly for small samples, of some of the most used
Markov chains order estimators.

It is quite intuitive that  the random information regarding the order of a Markov chain, is spread over
an exponentially growing set of empirical distributions $\Theta$ with $\vert \Theta \vert = m^{B+1}$,
where $\textbf{B}$ is the maximum integer $\eta$, as in $\alpha = (i_1i_2...i_\eta)$. It seems  reasonable to think that a 
small \textit{viable} sample, i.e. samples able to retrieve enough information to estimate the chain order, should have  size 
$\,\,n \approx O(m^{B+1}).$ Keeping in mind that for the present numerical simulation, the maximum length to 
be used is $B=5$, from now on the sample
sizes for $\vert E \vert = 3$ and $\vert E \vert = 4$ should be $n \approx 1.500$ and $n \approx 5.000$, respectively.

\vspace{20pt}
The following numerical simulation, based on an algorithm due to Raftery\cite{Raftery}, starts on with the generation 
of a  Markov chain transition matrix,  $\textbf{Q}=(q_{i_1i_2...i_{\kappa};i_{\kappa +1}})$ with entries 

\begin{linenomath*}
\begin{equation}
q_{i_1i_2...i_{\kappa};i_{\kappa +1}} = \sum_{t=1}^\kappa \lambda_{i_t} R(i_{\kappa +1},i_t), \,\, 
1 \leq i_t,i_{\kappa + 1} \leq m.
\label{markov_matrix_generation}
\end{equation}
\end{linenomath*}

where the matrix
\begin{linenomath*}
\begin{equation*}
{R}(i,j), \,\, 0 \leq i,j \leq m, \quad \sum_{i=1}^m R(i,j)=1, \,\, 1 \leq j \leq m
\end{equation*}
\end{linenomath*}

and the positive numbers 
\begin{linenomath*}
\begin{equation*}
\{ \lambda_i \}_{i=1}^\kappa, \,\, \sum_{i=1}^\kappa \lambda_i =1
\end{equation*}
\end{linenomath*}
are arbitrarily chosen in advance.

Once the  matrix $\textbf{Q}=(q_{i_1i_2...i_{\kappa};i_{\kappa +1}})$ is obtained, two hundreds
replications of the  Markov chain sample of size $n$, space state $E$ and transition matrix $\textbf{Q}$ 
are generated to  compare $GDL(\eta)$ performance against the standards, well known and 
already established order estimators just mentioned above.

Finally, after applying all estimators to
each one of the replicated samples, the final results two hundreds replications are registered in the form of tables.

\vspace{30pt}
{\textbf{Case I: \,\, Markov \,\, Chain \,\, Examples \,\, with }\,\,  $\kappa = 0, \,\, \vert E \vert = 3.$ }

\vspace{10pt}
Firstly, we choose the matrix $ \{Q_1, Q_2, Q_3\}$ to produce  samples with sizes $500 \leq n \leq 2.000$,
originated from Markov chains of order $\kappa=0$  with quite different probability distributions.

\begin{linenomath*}
\[ 
Q_{1}= 
\left[ 
\begin{array}{cccc}
\scriptstyle{0.33} & \scriptstyle{0.335} & \scriptstyle{0.335} \\
\scriptstyle{0.33} & \scriptstyle{0.335} & \scriptstyle{0.335} \\
\scriptstyle{0.33} & \scriptstyle{0.335} & \scriptstyle{0.335} 
\end{array} 
\right]
,\,
Q_{2}= 
\left[ 
\begin{array}{cccc}
\scriptstyle{0.05} & \scriptstyle{0.475} & \scriptstyle{0.475} \\
\scriptstyle{0.05} & \scriptstyle{0.475} & \scriptstyle{0.475} \\
\scriptstyle{0.05} & \scriptstyle{0.475} & \scriptstyle{0.475} 
\end{array} 
\right]
,\,
Q_{3}= 
\left[ 
\begin{array}{cccc}
\scriptstyle{0.05} & \scriptstyle{0.05} & \scriptstyle{0.90} \\
\scriptstyle{0.05} & \scriptstyle{0.05} & \scriptstyle{0.90} \\
\scriptstyle{0.05} & \scriptstyle{0.05} & \scriptstyle{0.90}  
\end{array} 
\right].
\]
\end{linenomath*}

\begin{tabular}{|c||c|c|c|c||c|c|c|c||c|c|c|c||c|c|c|c|}
                                                     \hline 
\multicolumn{13}{|c|}{\small 
     \textbf{\textit{ $\vert E \vert = 3 \qquad \leftrightarrow  \qquad \small  \kappa = 0  
                                         \qquad \leftrightarrow  \qquad \scriptstyle {\lambda_i\,=\,1/3, \,\, i\,=\,1,2,3.}$}}}\\
                                                     \hline
                                      &\multicolumn{4}{|c||}{\small $Q_1$}  
                                      &\multicolumn{4}{|c||}{\small $Q_1$}  
                                      &\multicolumn{4}{|c||}{\small $Q_1$}\\
                                                    \hline
                                      &\multicolumn{4}{|c||}{$\scriptstyle  n\,=\,500 $}
                                      &\multicolumn{4}{|c||}{$\scriptstyle  n\,=\,1.000 $}
                                      &\multicolumn{4}{|c||}{$\scriptstyle  n\,=\,1.500$ }\\
                                                    \hline \hline
$\scriptstyle k$     & $\scriptstyle Aic$ &  $\scriptstyle Bic$ & $\scriptstyle Edc$  & $\scriptstyle Gdl$ &
                       $\scriptstyle Aic$ &  $\scriptstyle Bic$ & $\scriptstyle Edc$  & $\scriptstyle Gdl$ &
                       $\scriptstyle Aic$ &  $\scriptstyle Bic$ & $\scriptstyle Edc$  & $\scriptstyle Gdl$ \\
                       	                             \hline \hline
$\scriptstyle 0$     & $\scriptstyle 75.5 \%$ & $\scriptstyle 100 \%$ & $\scriptstyle 100 \% $ & $\scriptstyle 99 \%  $ & 
                       $\scriptstyle 80   \%$ & $\scriptstyle 100 \%$ & $\scriptstyle 100 \% $ & $\scriptstyle 99.5 \%$ & 
                       $\scriptstyle 71.5 \%$ & $\scriptstyle 100 \%$ & $\scriptstyle 100 \% $ & $\scriptstyle 99 \%  $ \\
                                                            \hline 
$\scriptstyle 1$     & $\scriptstyle 24.5\% $ & $\scriptstyle     $ & $\scriptstyle     $ & $\scriptstyle 1 \%   $ & 
                       $\scriptstyle 18  \% $ & $\scriptstyle     $ & $\scriptstyle     $ & $\scriptstyle 0.5 \% $ & 
                       $\scriptstyle 22.5 \%$ & $\scriptstyle     $ & $\scriptstyle     $ & $\scriptstyle 1 \%   $ \\
                                                            \hline  
$\scriptstyle 2$     & $\scriptstyle        $ & $\scriptstyle  $ & $\scriptstyle    $ & $\scriptstyle    $ & 
                       $\scriptstyle  2 \%  $ & $\scriptstyle  $ & $\scriptstyle    $ & $\scriptstyle    $ & 
                       $\scriptstyle  6 \%  $ & $\scriptstyle  $ & $\scriptstyle    $ & $\scriptstyle    $ \\
                                                            \hline
$\scriptstyle 3$     & $\scriptstyle   $ & $\scriptstyle      $ & $\scriptstyle   $ & $\scriptstyle   $ & 
                       $\scriptstyle   $ & $\scriptstyle      $ & $\scriptstyle   $ & $\scriptstyle   $ & 
                       $\scriptstyle   $ & $\scriptstyle      $ & $\scriptstyle   $ & $\scriptstyle   $ \\
                                                            \hline
$\scriptstyle 4$     & $\scriptstyle \phantom{100\%}    $ & $\scriptstyle    $ & $\scriptstyle    $ & $\scriptstyle \phantom{100\%}    $ & 
                       $\scriptstyle \phantom{100\%}    $ & $\scriptstyle    $ & $\scriptstyle    $ & $\scriptstyle \phantom{100\%}    $ & 
                       $\scriptstyle \phantom{100\%}    $ & $\scriptstyle    $ & $\scriptstyle    $ & $\scriptstyle \phantom{100\%}    $ \\
                                                       \hline \hline
\end{tabular}

\vspace{10pt}
\begin{tabular}{|c||c|c|c|c||c|c|c|c||c|c|c|c||c|c|c|c|}
                                                     \hline 
\multicolumn{13}{|c|}{\small 
     \textbf{\textit{ $\vert E \vert = 3 \qquad \leftrightarrow  \qquad \small  \kappa = 0  
                                         \qquad \leftrightarrow  \qquad \scriptstyle {\lambda_i\,=\,1/3, \,\, i\,=\,1,2,3.}$}}}\\
                                                     \hline
                                      &\multicolumn{4}{|c||}{\small $Q_2$}  
                                      &\multicolumn{4}{|c||}{\small $Q_2$}  
                                      &\multicolumn{4}{|c||}{\small $Q_2$}\\
                                                     \hline
                                      &\multicolumn{4}{|c||}{$\scriptstyle  n\,=\,1.000 $}
                                      &\multicolumn{4}{|c||}{$\scriptstyle  n\,=\,1.500 $}
                                      &\multicolumn{4}{|c||}{$\scriptstyle  n\,=\,500   $}\\
                                                    \hline \hline
$\scriptstyle k$     & $\scriptstyle Aic$ &  $\scriptstyle Bic$ & $\scriptstyle Edc$  & $\scriptstyle Gdl$ &
                       $\scriptstyle Aic$ &  $\scriptstyle Bic$ & $\scriptstyle Edc$  & $\scriptstyle Gdl$ &
                       $\scriptstyle Aic$ &  $\scriptstyle Bic$ & $\scriptstyle Edc$  & $\scriptstyle Gdl$ \\
                       	                             \hline \hline
$\scriptstyle 0$     & $\scriptstyle 63.5\%$ & $\scriptstyle 100\% $ & $\scriptstyle 100\% $ & $\scriptstyle  99\% $ & 
                       $\scriptstyle 63\%  $ & $\scriptstyle 100\% $ & $\scriptstyle 100\% $ & $\scriptstyle  99\%   $ & 
                       $\scriptstyle 59\%  $ & $\scriptstyle 100\% $ & $\scriptstyle 100\% $ & $\scriptstyle  99\% $ \\
                                                            \hline
$\scriptstyle 1$     & $\scriptstyle 29\%  $ & $\scriptstyle     $ & $\scriptstyle     $ & $\scriptstyle  1\%  $ & 
                       $\scriptstyle 34.5\%$ & $\scriptstyle     $ & $\scriptstyle     $ & $\scriptstyle  1\%  $ & 
                       $\scriptstyle 37\%  $ & $\scriptstyle     $ & $\scriptstyle     $ & $\scriptstyle  1\%  $\\
                                                            \hline
$\scriptstyle 2$     & $\scriptstyle 7.5\% $ & $\scriptstyle        $ & $\scriptstyle       $ & $\scriptstyle    $ & 
                       $\scriptstyle 2.5\% $ & $\scriptstyle        $ & $\scriptstyle       $ & $\scriptstyle    $ & 
                       $\scriptstyle 4\%   $ & $\scriptstyle        $ & $\scriptstyle       $ & $\scriptstyle    $ \\
                                                            \hline
$\scriptstyle 3$     & $\scriptstyle      $ & $\scriptstyle       $ & $\scriptstyle       $ & $\scriptstyle      $ & 
                       $\scriptstyle      $ & $\scriptstyle       $ & $\scriptstyle       $ & $\scriptstyle      $ & 
                       $\scriptstyle      $ & $\scriptstyle       $ & $\scriptstyle       $ & $\scriptstyle      $ \\
                                                            \hline
$\scriptstyle 4$     & $\scriptstyle \phantom{100\%}    $ & $\scriptstyle    $ & $\scriptstyle    $ & $\scriptstyle \phantom{100\%}    $ & 
                       $\scriptstyle \phantom{100\%}    $ & $\scriptstyle    $ & $\scriptstyle    $ & $\scriptstyle \phantom{100\%}    $ & 
                       $\scriptstyle \phantom{100\%}    $ & $\scriptstyle    $ & $\scriptstyle    $ & $\scriptstyle \phantom{100\%}    $ \\
                                                       \hline \hline
\end{tabular}

\vspace{10pt}
\begin{tabular}{|c||c|c|c|c||c|c|c|c||c|c|c|c||c|c|c|c|}
                \hline 
\multicolumn{13}{|c|}{\small 
     \textbf{\textit{ $\vert E \vert = 3 \qquad \leftrightarrow  \qquad \small  \kappa = 0  
                                         \qquad \leftrightarrow  \qquad \scriptstyle {\lambda_i\,=\,1/3, \,\, i\,=\,1,2,3.}$}}}\\
                                                     \hline
                                      &\multicolumn{4}{|c||}{\small $Q_3$}  
                                      &\multicolumn{4}{|c||}{\small $Q_3$}  
                                      &\multicolumn{4}{|c||}{\small $Q_3$}\\
                                                     \hline
                                      &\multicolumn{4}{|c||}{$\scriptstyle  n\,=\,1.000 $}
                                      &\multicolumn{4}{|c||}{$\scriptstyle  n\,=\,1.500 $}
                                      &\multicolumn{4}{|c||}{$\scriptstyle  n\,=\,2.000$ }\\
                                                    \hline \hline
$\scriptstyle k$     & $\scriptstyle Aic$ &  $\scriptstyle Bic$ & $\scriptstyle Edc$  & $\scriptstyle Gdl$ &
                       $\scriptstyle Aic$ &  $\scriptstyle Bic$ & $\scriptstyle Edc$  & $\scriptstyle Gdl$ &
                       $\scriptstyle Aic$ &  $\scriptstyle Bic$ & $\scriptstyle Edc$  & $\scriptstyle Gdl$ \\
                       	                             \hline \hline
$\scriptstyle 0$     & $\scriptstyle 43\%$ & $\scriptstyle 100\%$ & $\scriptstyle 100\% $ & $\scriptstyle 98\%$ & 
                       $\scriptstyle 47\%$ & $\scriptstyle 100\%$ & $\scriptstyle 99.5\%$ & $\scriptstyle 96\%$ & 
                       $\scriptstyle 46\%$ & $\scriptstyle 100\%$ & $\scriptstyle 100\% $ & $\scriptstyle 97\%$ \\
                                                            \hline
$\scriptstyle 1$     & $\scriptstyle 53\%  $ & $\scriptstyle     $ & $\scriptstyle      $ & $\scriptstyle 2\%  $ & 
                       $\scriptstyle 51.5\%$ & $\scriptstyle     $ & $\scriptstyle 0.5\%$ & $\scriptstyle 4\%  $ & 
                       $\scriptstyle 50.5\%$ & $\scriptstyle     $ & $\scriptstyle      $ & $\scriptstyle 2\%  $ \\
                                                            \hline
$\scriptstyle 2$     & $\scriptstyle  4\% $ & $\scriptstyle     $ & $\scriptstyle    $ & $\scriptstyle      $ & 
                       $\scriptstyle 1.5\%$ & $\scriptstyle     $ & $\scriptstyle    $ & $\scriptstyle      $ & 
                       $\scriptstyle 3.5\%$ & $\scriptstyle     $ & $\scriptstyle    $ & $\scriptstyle  1\% $ \\
                                                            \hline
$\scriptstyle 3$     & $\scriptstyle      $ & $\scriptstyle       $ & $\scriptstyle     $ & $\scriptstyle        $ & 
                       $\scriptstyle      $ & $\scriptstyle       $ & $\scriptstyle     $ & $\scriptstyle        $ & 
                       $\scriptstyle      $ & $\scriptstyle       $ & $\scriptstyle     $ & $\scriptstyle        $ \\
                                                            \hline
$\scriptstyle 4$     & $\scriptstyle \phantom{100\%} $ & $\scriptstyle    $ & $\scriptstyle    $ & $\scriptstyle \phantom{100\%}    $ & 
                       $\scriptstyle \phantom{100\%} $ & $\scriptstyle    $ & $\scriptstyle    $ & $\scriptstyle \phantom{100\%}    $ & 
                       $\scriptstyle \phantom{100\%} $ & $\scriptstyle    $ & $\scriptstyle    $ & $\scriptstyle \phantom{100\%}    $ \\
                                                       \hline \hline
\end{tabular}

\vspace{25pt}
Notice that for a fixed sample size $n = \{ 500, 1.000, 1.500, 2.000 \}$, the order estimator $\widehat{\kappa}_{AIC}$ 
steadily overestimate the real order $\kappa = 0$ with the excessiveness depending on the  probability distribution of the Markov chain. 
Differently, the order estimators $\widehat{\kappa}_{BIC}$, $\widehat{\kappa}_{EDC}$ and $\widehat{\kappa}_{GDL}$ show consistent performance, 
mainly obtaining the right order, free from the influence of the sample size and the generating matrix.
Regarding $\widehat{\kappa}_{BIC}$ and  $\widehat{\kappa}_{EDC}$ improved effect, most likely depends on their correcting factor,
$\frac{\log(n)}{2}$ and $\left( \frac{\log\log(n)}{2 (\vert E \vert - 1)} \right )$
which tend to decrease the estimated order.

For $\vert E \vert = 4$ the greater complexity of a Markov chain of order $\kappa = 3$ impose the use of larger sample size
for estimators to acomplish some reliability.
Finally, we choose the matrix $ \{ Q_6, Q_7 \}$ to produce samples with size $n = 5.000$,
originated from Markov chains of order $\kappa \in \{ 2,3,0 \}$  like in the previous cases.

\vspace{30pt}
\begin{linenomath*}
\[ 
Q_{6}= 
\left[ 
\begin{array}{cccc}
\scriptstyle{0.05} & \scriptstyle{0.05} & \scriptstyle{0.05} & \scriptstyle{0.85} \\
\scriptstyle{0.05} & \scriptstyle{0.05} & \scriptstyle{0.85} & \scriptstyle{0.05} \\
\scriptstyle{0.05} & \scriptstyle{0.85} & \scriptstyle{0.05} & \scriptstyle{0.05} \\
\scriptstyle{0.85} & \scriptstyle{0.05} & \scriptstyle{0.05} & \scriptstyle{0.05}
\end{array} 
\right]
,\qquad
Q_{7}= 
\left[ 
\begin{array}{cccc}
\scriptstyle{0.05} & \scriptstyle{0.05} & \scriptstyle{0.05} & \scriptstyle{0.85} \\
\scriptstyle{0.05} & \scriptstyle{0.05} & \scriptstyle{0.05} & \scriptstyle{0.85} \\
\scriptstyle{0.05} & \scriptstyle{0.05} & \scriptstyle{0.05} & \scriptstyle{0.85} \\
\scriptstyle{0.05} & \scriptstyle{0.05} & \scriptstyle{0.05} & \scriptstyle{0.85}
\end{array} 
\right].
\]
\end{linenomath*}

\begin{tabular}{|c||c|c|c|c||c|c|c|c||c|c|c|c||c|c|c|c|}
                \hline 
\multicolumn{13}{|c|}{\small \textbf{\textit{ $\vert E \vert \,=\, 4 \,\,\,\,\,\,\, \leftrightarrow \,\,\,\,\,\,\,\, n\,=\, 5.000   $}} } \\
                                                     \hline
     &\multicolumn{4}{|c||}{\small $Q_6  \Leftrightarrow \,\, \scriptstyle{ \lambda_i\,=\,1/2, \,\, i\,=\,1,2.  }$}
     &\multicolumn{4}{|c||}{\small $Q_6  \Leftrightarrow \,\, \scriptstyle{ \lambda_i\,=\,1/3, \,\, i\,=\,1,2,3.}$}
     &\multicolumn{4}{|c|}{\small  $Q_7  \Leftrightarrow \,\, \scriptstyle{ \lambda_i\,=\,1/3, \,\, i\,=\,1,2,3.}$}\\
                                                      \hline
                                      &\multicolumn{4}{|c||}{\small $ \kappa\,=\, 2 $}
                                      &\multicolumn{4}{|c||}{\small $ \kappa\,=\, 3 $}
                                      &\multicolumn{4}{|c|}{\small  $ \kappa\,=\, 0 $}\\
                                                    \hline \hline
$\scriptstyle k$     & $\scriptstyle Aic$ &  $\scriptstyle Bic$ & $\scriptstyle Edc$  & $\scriptstyle Gdl$ &
                       $\scriptstyle Aic$ &  $\scriptstyle Bic$ & $\scriptstyle Edc$  & $\scriptstyle Gdl$ &
                       $\scriptstyle Aic$ &  $\scriptstyle Bic$ & $\scriptstyle Edc$  & $\scriptstyle Gdl$ \\
                       	                             \hline \hline
$\scriptstyle 0$     & $\scriptstyle    $ & $\scriptstyle    $ & $\scriptstyle         $ & $\scriptstyle      $ & 
                       $\scriptstyle    $ & $\scriptstyle    $ & $\scriptstyle         $ & $\scriptstyle      $ & 
                       $\scriptstyle 85\%$ & $\scriptstyle 100\%$ & $\scriptstyle 100\%$ & $\scriptstyle 100\%$ \\
                                                            \hline
$\scriptstyle 1$     & $\scriptstyle     $ & $\scriptstyle    $ & $\scriptstyle    $ & $\scriptstyle    $ & 
                       $\scriptstyle     $ & $\scriptstyle    $ & $\scriptstyle    $ & $\scriptstyle    $ & 
                       $\scriptstyle 15\%$ & $\scriptstyle    $ & $\scriptstyle    $ & $\scriptstyle    $ \\
                                                            \hline
$\scriptstyle 2$     & $\scriptstyle 100\%$ & $\scriptstyle 100\% $ & $\scriptstyle 100\%$ & $\scriptstyle 100\%$ & 
                       $\scriptstyle      $ & $\scriptstyle   99\%$ & $\scriptstyle      $ & $\scriptstyle   4\%$ & 
                       $\scriptstyle      $ & $\scriptstyle       $ & $\scriptstyle      $ & $\scriptstyle      $ \\
                                                            \hline
$\scriptstyle 3$     & $\scriptstyle      $ & $\scriptstyle      $ & $\scriptstyle      $ & $\scriptstyle       $ & 
                       $\scriptstyle 100\%$ & $\scriptstyle   1\%$ & $\scriptstyle 100\%$ & $\scriptstyle   96\%$ & 
                       $\scriptstyle      $ & $\scriptstyle      $ & $\scriptstyle      $ & $\scriptstyle       $ \\
                                                            \hline
$\scriptstyle 4$     & $\scriptstyle    $ & $\scriptstyle    $ & $\scriptstyle    $ & $\scriptstyle    $ & 
                       $\scriptstyle    $ & $\scriptstyle    $ & $\scriptstyle    $ & $\scriptstyle    $ & 
                       $\scriptstyle    $ & $\scriptstyle    $ & $\scriptstyle    $ & $\scriptstyle    $ \\
                                                            \hline
$\scriptstyle 5$     & $\scriptstyle    $ & $\scriptstyle    $ & $\scriptstyle    $ & $\scriptstyle    $ & 
                       $\scriptstyle    $ & $\scriptstyle    $ & $\scriptstyle    $ & $\scriptstyle    $ & 
                       $\scriptstyle    $ & $\scriptstyle    $ & $\scriptstyle    $ & $\scriptstyle    $ \\
                                                            \hline
$\scriptstyle 6$     & $\scriptstyle    $ & $\scriptstyle    $ & $\scriptstyle    $ & $\scriptstyle    $ & 
                       $\scriptstyle    $ & $\scriptstyle    $ & $\scriptstyle    $ & $\scriptstyle    $ & 
                       $\scriptstyle    $ & $\scriptstyle    $ & $\scriptstyle    $ & $\scriptstyle    $ \\ 
                                                            \hline\hline
\end{tabular}

\vspace{25pt}
For the order for $\vert E \vert = 4$, $\kappa =0$, apparently  $\widehat{\kappa}_{AIC}$ keeps overestimating 
the order in some degree, while $\widehat{\kappa}_{BIC}$ as in example $\kappa = 3$
severely underestimate the order, presumably due to the excessive weight of the correcting factors $\frac{\log(n)}{2}$. 
On the contrary $\widehat{\kappa}_{EDC}$ and $\widehat{\kappa}_{GDL}$ behaves quite well in same setting.

\end{document}